\documentclass[12pt]{article} \usepackage{amsfonts}
\begin{document}
\baselineskip=14pt
\pagestyle{plain}
{\Large
\newcommand{\mun}{\mu_{n,j}}
\newcommand{\aln}{\alpha_{n,j}(\mu)}
\newcommand{\wnj}{W_{N,j}(\mu)}
\newcommand{\ntm}{|2n-\theta-\mu|}
\newcommand{\pwp}{PW_\pi^-}
\newcommand{\cpm}{\cos\pi\mu}
\newcommand{\spm}{\frac{\sin\pi\mu}{\mu}}
\newcommand{\tet}{(-1)^{\theta+1}}
\newcommand{\dm}{\Delta(\mu)}
\newcommand{\dnm}{\Delta_N(\mu)}
\newcommand{\sni}{\sum_{n=N+1}^\infty}
\newcommand{\muk}{\sqrt{\mu^2+q_0}}
\newcommand{\agt}{\alpha, \gamma, \theta,}
\newcommand{\muq}{\sqrt{\mu^2-q_0}}
\newcommand{\smn}{\sum_{n=1}^\infty}
\newcommand{\lop}{{L_2(0,2\pi)}}
\newcommand{\tdm}{\tilde\Delta_+(\mu_n)}
\centerline{\bf Inverse problems of spectral analysis for the Sturm-Liouville}
\centerline {\bf operator with regular boundary conditions}

\medskip
\medskip
\centerline { Alexander Makin}
\medskip
\medskip

\centerline{Abstract. We consider the Sturm-Liouville operator}
\centerline{$Lu=u''-q(x)u$ defined on $(0,\pi)$ with
regular}\centerline{ but not strongly regular boundary conditions.
Under}\centerline{ some supplementary assumptions we prove that
the}\centerline{set of potentials $q(x)$ that ensure an
asymptotically}\centerline{ multiple spectrum is everywhere dense
in $L_1(0,\pi)$.}

\medskip
\medskip

In the present paper we consider eigenvalue problems for the Sturm-Liouville equation
$$
u''-q(x)u+\lambda u=0 \eqno (1)
$$
with two-point boundary conditions
$$
B_i(u)=a_{i1}u'(0)+a_{i2}u'(\pi)+a_{i3}u(0)+a_{i4}u(\pi)=0, \eqno (2)
$$
where $B_i(u)$ $(i=1,2)$ are linearly independent forms with arbitrary complex-valued
coefficients. Function $q(x)$ is an arbitrary complex-valued function of
class $L_2(0,\pi)$.

It is convenient to write conditions (2) in the matrix form
$$
A=\left(
\begin{array}{cccc}
a_{11}&a_{12}&a_{13}&a_{14}\\
a_{21}&a_{22}&a_{23}&a_{24}
\end{array}
\right)
$$
and denote the matrix composed of the $i$th and $j$th columns of $A$
 $(1\le i<j\le4)$ by $A(ij)$; we set $A_{ij}=det$ $A(ij)$.
Let the boundary conditions (2) be regular but not strongly
regular [1, pp. 71-73], which, by [1, p. 73] is equivalent to the
conditions
$$
A_{12}=0, \quad A_{14}+A_{23}\ne0,\quad A_{14}+A_{23}=\mp(A_{13}+A_{24}).\eqno(3)
$$

To investigate this class of problems it is appropriate [2] to divide
conditions (2) satisfying (3) into 4 types:

I) $A_{14}=A_{23}$, $A_{34}=0$;

II) $A_{14}=A_{23}$, $A_{34}\ne 0$;

III) $A_{14}\ne A_{23}$, $A_{34}=0$;

IV) $A_{14}\ne A_{23}$, $A_{34}\ne 0$.

An eigenvalue problem for equation (1) with boundary conditions of
type I, II, III or IV is called a problem of type I, II, III, IV,
respectively. It is shown [2] that any boundary conditions of type
I are equivalent to the boundary conditions specified by matrix
$A$, where
$$
A=\left(
\begin{array}{cccc}
1&-1&0&0\\
0&0&1&-1
\end{array}
\right) \quad or\quad
A=\left(
\begin{array}{cccc}
1&1&0&0\\
0&0&1&1
\end{array}
\right),
$$
i.e. to periodic or antiperiodic boundary conditions.
Inverse problems for this case have been studied well (see, for instance, [3],
 [4]).

It also shown in [2] that any boundary conditions of type II are equivalent
to the boundary conditions specified by matrix
 $A$, where
$$
A=\left(
\begin{array}{cccc}
1&-1&0&a_{14}\\
0&0&1&-1
\end{array}
\right) \quad or\quad
A=\left(
\begin{array}{cccc}
1&1&0&a_{14}\\
0&0&1&1
\end{array}
\right),
$$
and in both cases $a_{14}\ne0$. If $a_{14}$ is a real number and
$q(x)$ is a real-valued function, then the corresponding boundary
value problem is selfadjoint [2]. For this case inverse problem
was studied in [5], [6].

Exhausting description of boundary conditions of type III and IV
is given in [2]. In particular, it is known that  all of them are
nonselfadjoint. The main purpose of this paper is to investigate
inverse problems generated by boundary conditions of type III or
IV.

Denote by $c(x,\mu), s(x,\mu)$ $(\lambda=\mu^2)$ a fundamental
system of solutions to the equation (1) with the boundary
conditions $c(0,\mu)=s'(0,\mu)=1$, $c'(0,\mu)=s(0,\mu)=0$. The
following identity is well known
$$
c(x,\mu)s'(x,\mu)-c'(x,\mu)s(x,\mu)=1.\eqno(4)
$$
Suppose condition (3) holds. Then simple computations show that
the characteristic equation of any problem (1)+(2)  can be reduced
to the form $\Delta(\mu)=0$, where
$$
\Delta(\mu)=(-1)^{\theta+1}+\alpha c(\pi,\mu)+(1-\alpha)s'(\pi,\mu)+\gamma s(\pi,\mu), \eqno(5)
$$
where $\alpha=A_{14}/(A_{14}+A_{23})$,
$\gamma=-A_{34}/(A_{14}+A_{23})$, and $\theta=0$ if the sign "-"
is situated in (3) (case 1), and $\theta=1$ if the sign
 "+" is situated in (3) (case 2). Thus, for boundary conditions of type I $\alpha=1/2, \gamma=0$,
for ones of type II $\alpha=1/2, \gamma\ne0$, for ones of type III
$\alpha\ne1/2, \gamma=0$, for ones of type IV $\alpha\ne1/2,
\gamma\ne0$. Obviously, the spectrum of problem (1)+(2) is defined
uniquely by a quadruple of parameters $(\alpha, \gamma, \theta,
q(x))$. The characteristic determinant $\Delta(\mu)$ of problem
(1)+(2) defined by formula (5) is called the characteristic
determinant corresponding to a collection $(\alpha, \gamma,
\theta, q(x))$. Throughout the following $||f||$ stands for
$||f||_{L_2(0,\pi)}$, $<q>=\frac{1}{\pi}\int_0^\pi q(x)dx$.
 By $\Gamma(z,r)$ we denote a disk of radius $r$ centered at $z$.
By $PW_\sigma$ we denote the class of all entire functions $f(z)$ of exponential type
not exceeding $\sigma$ such that $||f(z)||_{L_2(R)}<\infty$, and by $PW_\sigma^-$
we denote the set of odd functions belonging to $PW_\sigma$.

The following theorem establishes the necessary  conditions for a
characteristic determinant $\Delta(\mu)$.

{\bf Theorem 1.} {\it If a function $\Delta(\mu)$ is a
characteristic determinant corresponding to a collection $(\alpha,
\gamma, \theta, q(x))$ then
$$
\Delta(\mu)=(-1)^{\theta+1}+\cos \pi\mu+(\gamma+\frac{\pi<q>}{2})\frac{\sin\pi\mu}{\mu}+\frac{f(\mu)}{\mu},
$$
where $f(\mu)\in PW_\pi^-$.}

Proof. Let $e(x,\mu)$ be the solution to equation (1) satisfying
initial conditions $e(0,\mu)=1$, $e'(0,\mu)=i\mu$ and let
$K(x,t)$, $K^+(x,t)=K(x,t)+K(x,-t)$, $K^-(x,t)=K(x,t)-K(x,-t)$ be
the kernels of transformation [3, pp. 17-18] realizing the
representations
$$
e(x,\mu)=e^{i\mu x}+\int_{-x}^xK(x,t)e^{i\mu t}dt,
$$
$$
c(x,\mu)=\cos\mu x+\int_0^xK^+(x,t)\cos\mu tdt,
$$
$$
s(x,\mu)=\frac{\sin\mu x}{\mu}+\int_0^xK^-(x,t)\frac{\sin\mu t}{\mu}dt.\eqno(6)
$$
It was shown in [7], [8] that
$$
c(\pi,\mu)=\cos\pi\mu+\frac{\pi}{2}<q>\spm-\int_0^\pi\frac{\partial K^+(\pi,t)}{\partial t}\frac{\sin\mu t}{\mu}dt,\eqno(7)
$$
$$
s(\pi,\mu)=\frac{\sin\pi\mu}{\mu}-\frac{\pi}{2\mu^2}<q>\cpm+\int_0^\pi\frac{\partial
K^-(\pi,t)}{\partial t}\frac{\cos\mu t}{\mu^2}dt,\eqno(8)
$$

$$
s'(\pi,\mu)=\cos\pi\mu+\frac{\pi}{2}<q>\spm+\int_0^\pi\frac{\partial K^-(\pi,t)}{\partial x}\frac{\sin\mu t}{\mu}dt.\eqno(9)
$$
Substituting the right-hand parts of (6), (7), (9) in (5), we
obtain
$$
\begin{array}{c}
\dm=\tet+\cpm+(\gamma+\frac{\pi}{2}<q>)\spm+\\
+\frac{1}{\mu}\int_0^\pi[-\alpha\frac{\partial K^+(\pi,t)}{\partial t}+(1-\alpha)\frac{\partial K^-(\pi,t)}{\partial x}+\gamma K^-(x,t)]\sin\mu tdt.
\end{array}
$$
It follows from the last equality and the Paley-Wiener theorem that theorem 1 is valid.

{\bf Theorem 2.} {\it Let a function $u(\mu)$ have the form
$$
u(\mu)=(-1)^{\theta+1}+\cos \pi\mu+(\gamma+\frac{\pi q_0}{2})\frac{\sin\pi\mu}{\mu}+\frac{f(\mu)}{\mu},\eqno(10)
$$
where $f(\mu)\in PW_\pi^-$. Let $\gamma$, $q_0$, $\alpha$ be
arbitrary complex numbers but $\alpha\ne1/2$, $\alpha\ne0$,
$\alpha\ne1$. Then there exists a function $q(x)\in L_2(0,\pi)$
such that corresponding to the collection $(\alpha, \gamma,
\theta, q(x))$ the characteristic determinant $\Delta(\mu)=u(\mu)$
and $<q>=q_0$.}

Proof. At first we consider the case $q_0=0$. Denote
$u_+(\mu)=u(\mu)-\tet$. Let $\varepsilon_1$ be an arbitrary
positive number. Since [9, pp. 115, 125]
$$|f(\mu)|\le C_1||f(\mu)||_{L_2(R)}e^{\pi|Im\mu|},\eqno(11)$$
we see that there exists a number $N$ large enough for the
inequality

$$
|u_+(\mu)-\cpm|<\varepsilon_1\eqno(12)
$$
to be valid on the set $|Im\mu|\le1$, $Re\mu\ge N$.

Let $\{\mu_n\}$ $(n=1,2,\ldots)$ be a strictly monotone increasing
sequence of positive numbers such that
$|\mu_n-(N+1/2)|<\varepsilon_1$ if $1\le n\le N$ and $\mu_n=n$ if
$n\ge N+1$. Similarly [10] let us consider the function
$$
s(\mu)=\pi\prod_{n=1}^\infty\frac{\mu_n^2-\mu^2}{n^2}=\spm\prod_{n=1}^N\frac{\mu_n^2-\mu^2}{n^2-\mu^2}.\eqno(13)
$$
It is evident that all zeros of $s(\mu)$ are simple and for any
$n$ we have the inequality
$$
(-1)^n\dot s(\mu_n)>0.\eqno(14)
$$
It is shown in [10] that
$$
\dot s(n)=\frac{\pi(-1)^n}{n}(1+C_0n^{-2}+O(n^{-4})),\eqno(15)
$$
where $C_0$ is a constant. It follows from [10] that

$$
s(\mu)=\spm+O(\mu^{-3})\eqno(16)
$$
if $|Im\mu|\le1$.

 Let us consider the equation
$$
\alpha z^2-u_+(\mu_n)z+(1-\alpha)=0.\eqno(17)
$$
It has the roots
$$
c_n^\pm=\frac{u_+(\mu_n)\pm\sqrt{u_+^2(\mu_n)-4\alpha(1-\alpha)}}{2\alpha}.\eqno(18)
$$
We denote $\tilde c^\pm=\pm\sqrt{-4\alpha(1-\alpha)}/(2\alpha)$. Obviously,

$\tilde c^+=-\tilde c^-\ne0$. Since the functions
$\frac{z\pm\sqrt{z^2-4\alpha(1-\alpha)}}{2\alpha}$ are continuous
in the neighborhood of zero, we see that for any $\varepsilon>0$
there exists $\sigma>0$ such that for any $|z|<\sigma$
$$
|\frac{z+\sqrt{z^2-4\alpha(1-\alpha)}}{2\alpha}-\tilde c^+|<\varepsilon,
$$
$$
|\frac{z-\sqrt{z^2-4\alpha(1-\alpha)}}{2\alpha}-\tilde c^-|<\varepsilon.
$$
Since $\alpha\ne1/2$, we see that the functions
$$\begin{array}{c}
g_n^\pm(z)=\frac{(-1)^n+z\pm\sqrt{[(-1)^n+z]^2-4\alpha(1-\alpha)}}{2\alpha}=\\
=\frac{(-1)^n+z\pm\sqrt{(1-2\alpha)^2+2(-1)^nz+z^2}}{2\alpha}\end{array}
$$
are continuous in the neighborhood of zero. For definiteness we
will count that $\sqrt{(1-2\alpha)^2}=1-2\alpha$. This yields that
for any $\varepsilon>0$ there exists $\delta>0$ such that for any
$|z|<\delta$ $|g_n^-(z)-1|<\varepsilon$ if $n$ is even, and
$|g_n^+(z)+1|<\varepsilon$ if $n$ is odd. Let
$\varepsilon_1<\delta$ and $\varepsilon_1(1+\pi)<\sigma$. It
follows from (12) that for $n=1,\ldots,N$ all the numbers $c_n^+$
are contained in the disk $\Gamma(\tilde c^+,\varepsilon)$, and
all the numbers $c_n^-$ are contained in the disk $\Gamma(\tilde
c^-,\varepsilon)$.

If $n\ge N+1$, then for even $n$ all the numbers $c_n^-$ are
contained in the disk $\Gamma(1, \varepsilon)$, and for odd $n$
all the numbers $c_n^+$ are contained in the disk
$\Gamma(-1,\varepsilon)$. Let $c_n$ be a root of equation (17)
chosen according to the rule formulated below.

Let us consider two cases.

1) the numbers $c_n^+$ and $c_n^-$ do not lie on the imaginary
axis. Then we choose $\varepsilon$ $(\varepsilon<1/2)$ small
enough for one of the disks $\Gamma(\tilde c^+, \varepsilon)$,
$\Gamma(\tilde c^-,\varepsilon)$ to be lain strictly in the right
half-plane, and the other disk to be lain strictly in the left
half-plane. If $n=1,\ldots,N$, then for even $n$ we choose $c_n$
so as the point $c_n$ belongs to one of the disks mentioned above
which lies in the right half-plane, and for odd $n$ we choose
$c_n$ so as the point $c_n$ belongs to one of the disks mentioned
above which lies in the left half-plane. If $n\ge N+1$ we count
that for even $n$ $c_n=c_n^-$ and for odd $n$ $c_n=c_n^+$. Then
for any $n=1, 2,\ldots$ $(-1)^nRec_n>0$. It follows from this and
(14) that for any $n$ $Rez_n>0$, where
$$
z_n=\frac{c_n}{\mu_n\dot s(\mu_n)}.\eqno(19)
$$

2) the numbers $\tilde c^+$ and $\tilde c^-$ lie on the imaginary
axis. Let a line $l$ be a line which passes through the origin and
which is parallel to the line which passes through the points $1$
and $\tilde c^+$. Then these points lie in the same half-plane
relative to the line $l$, hence, a number $\varepsilon$ can be
chosen small enough for the disks $\Gamma(1,\varepsilon)$,
$\Gamma(\tilde c^+,\varepsilon)$ also to be lain strictly in the
same half-plane relative to the line $l$, therefore, by virtue of
symmetry, the disks $\Gamma(-1,\varepsilon)$, $\Gamma(\tilde
c_n^-,\varepsilon)$ lie strictly on the other half-plane relative
to the line $l$. If $n=1,\ldots,N$ then for even $n$ we choose
$c_n$ so as the point $c_n$ belongs to the disk $\Gamma(\tilde
c^+,\varepsilon)$, and for odd $n$ we choose $c_n$ so as the point
$c_n$ belongs to the disk $\Gamma(\tilde c^-,\varepsilon)$. If
$n\ge N+1$, then for even $n$ we count that $c_n=c_n^-$, and for
odd $n$ we count that $c_n=c_n^+$. Then for all even $n$ the
points $c_n$ lie strictly in the same half-plane relative to the
line $l$ and for all odd $n$ the points $c_n$ lie strictly in the
other half-plane relative to the line $l$. This, together with
(14), yields that for all $n$ the numbers $z_n$ lie strictly in
the same half-plane relative to the line $l$.

Let us set $F(x,t)=F_0(x,t)+\hat F(x,t)$, where
$$
F_0(x,t)=\sum_{n=1}^N\left(\frac{2c_n}{\mu_n\dot s(\mu_n)}
\sin\mu_nx\sin\mu_nt-
\frac{2}{\pi}\sin nx\sin nt\right),
$$
$$
\hat F(x,t)=\sum_{n=N+1}^\infty\left(\frac{2c_n}{\mu_n\dot s(\mu_n)}
\sin\mu_nx\sin\mu_nt-
\frac{2}{\pi}\sin nx\sin nt\right).\eqno(20)
$$
Clearly, $F_0(x,t)\in C^\infty(R^2)$. Let us consider the function
$\hat F(x,t)$. For convenience we denote
$R(n)=2(-1)^nf(n)/n+f^2(n)/n^2$. If $n\ge N+1$, then from (18) and
the rule of choice of the roots of equation(17) we obtain
$$
\begin{array}{c}
c_n=\frac{(-1)^n+f(n)/n-(-1)^n\sqrt{[(-1)^n+f(n)/n]^2-4\alpha(1-\alpha)}}{2\alpha}=\\
=\frac{(-1)^n+f(n)/n-(-1)^n\sqrt{(1-2\alpha)^2+R(n)}}{2\alpha}=\\
=\frac{(-1)^n+f(n)/n-(-1)^n[(1-2\alpha)+R(n)/(2(1-2\alpha))+O(R^2(n))]}{2\alpha}=\\
=(-1)^n-\frac{f(n)}{(1-2\alpha)n}+O(f^2(n)/n^2).
\end{array}\eqno(21)
$$
It follows from (11), (15), (20), (21) that
$$
\begin{array}{c}
\hat F(x,t)=\sum_{n=N+1}^\infty\frac{2}{\pi}\left(\frac{1-(-1)^n\frac{f(n)}{(1-2\alpha)n}+O(f^2(n)/n^2)}{1+c_0/n^2+O(1/n^4)}-1
\right)\sin nx\sin nt=\\
=\sum_{n=N+1}^\infty\frac{2}{\pi}[(1-(-1)^n\frac{f(n)}{(1-2\alpha)n}+O(1/n^2))\times\\(1-c_0/n^2+O(1/n^4))-1
]\sin nx\sin nt=\\
=\frac{2}{\pi}\sum_{n=N+1}^\infty((-1)^{n+1}\frac{f(n)}{(1-2\alpha)n}+O(1/n^2))\sin nx\sin nt=\\
=(\hat G(x-t)-\hat G(x+t))/2,

\end{array}
$$
where $$\hat
G(y)=\frac{2}{\pi}\sum_{n=N+1}^\infty((-1)^{n+1}\frac{f(n)}{(1-2\alpha)n}+O(1/n^2))\cos
ny.$$ By virtue of the Paley-Wiener theorem and the Parseval
equality, we have
$$
\sum_{n=1}^\infty|f(n)|^2=\frac{1}{2}||f(\mu)||_{L_2(R)},
$$ hence,
$\hat G(y)\in W_2^1[0,2\pi]$. Thus, we obtain
$$
F(x,t)=F_0(x,t)+(\hat G(x-t)-\hat G(x+t))/2,\eqno(22)
$$
where the functions $F_0(x,t)$ and $\hat G(y)$ belong to the classes mentioned above.

Let us now consider the Gelfand-Levitan equation
$$
K(x,t)+F(x,t)+\int_0^xK(x,t)F(s,t)ds=0\eqno(23)
$$
and prove that it is uniquely solvable in $L_2(0,x)$ for every $x\in[0,\pi]$.
To this end it is sufficient to prove that the corresponding homogeneous equation has
the trivial solution only.

Let $f(t)\in L_2(0,x)$ and
$$
f(t)+\int_0^xF(s,t)f(s)ds=0.
$$
Similarly [10] multiplying this equation by $\bar f(t)$ and
integrating the resulting equation over the $[0,x]$, we obtain
$$\begin{array}{c}
\int_0^x|f(t)|^2dt+
\sum_{n=1}^\infty\frac{2c_n}{\mu_n\dot s(\mu_n)}\int_0^x\bar f(t)\sin\mu_ntdt\int_0^xf(s)\sin\mu_nsds-\\
-\sum_{n=1}^\infty\frac{2}{\pi}\int_0^x\bar f(t)\sin ntdt\int_0^xf(s)\sin nsds=0.
\end{array}
$$
This, together with the Parseval equality for the system $\{\sin nt\}_1^\infty$
 on the segment $[0,\pi]$, yields
$$
\sum_{n=1}^\infty z_n|\int_0^xf(t)\sin\mu_ntdt|^2=0,
$$
where the numbers $z_n$ are defined by (19). Since all $z_n$ are
situated strictly in the same half-plane relative to a line which
passes through the origin, we see that
$\int_0^xf(t)\sin\mu_ntdt=0$ for any $n=1,2,\ldots$. Since [11,
12] the system $\{\sin\mu_nt\}_1^\infty$ is complete on the
segment $[0,\pi]$, it follows that $f(t)\equiv0$ on the segment
$[0,x]$.

Let $\hat K(x,t)$ be the unique solution of equation (23). Let us
set $\hat q(x)=2\frac{d}{dx}\hat K(x,x)$. It follows from (22)
[10] that $\hat q(x)\in L_2(0,\pi)$. We denote by $\hat s(x,\mu)$,
$\hat c(x,\mu)$ the fundamental system of solutions to equation
(1) with the potential $\hat q(x)$ and the initial conditions
$\hat s(0,\mu)=\hat c'(0,\mu)=0$, $\hat c(0,\mu)=\hat
s'(0,\mu)=1$. Repeating the reasons of [10], we obtain that $\hat
s(\pi,\mu)\equiv s(\mu)$, therefore, the numbers $\mu_n^2$ form
the spectrum of the Dirichlet problem for equation (1) with the
potential $\hat q(x)$. Similarly [10] we  also obtain that $\hat
c(\pi,\mu_n)=c_n$. This, together with identity (4), yields $\hat
s'(\pi,\mu_n)=1/c_n$.

It follows from (8), (16) and the Riemann lemma [3, p. 36] that
$<\hat q>=0$.

Let $\hat\Delta(\mu)$ be the characteristic determinant
corresponding to the collection $(\alpha,\gamma,\theta,\hat
q(x))$. Let us prove that $\hat\Delta(\mu)\equiv u(\mu)$.
According to theorem 1, the function $\hat \Delta(\mu)$ may be
represented in the form
$$
\hat\Delta(\mu)=(-1)^{\theta+1}+\cos \pi\mu+\gamma\frac{\sin\pi\mu}{\mu}+\frac{\hat f(\mu)}{\mu},
$$
where $\hat f(\mu)\in PW_\pi^-$. Taking into account (5) and
remembering that the numbers $c_n$ are the roots of equation (17),
we obtain
$$
\begin{array}{c}
\hat\Delta(\mu_n)=\tet+
\alpha \hat c(\pi,\mu_n)+(1-\alpha)\hat s'(\pi,\mu_n)+\gamma\hat s(\pi,\mu_n)=\\
=\tet+\alpha c_n+(1-\alpha)/c_n=\tet+u_+(\mu_n)=u(\mu_n).
\end{array}
$$
This implies that the function
$$
\Phi(\mu)=\frac{u(\mu)-\hat\Delta(\mu)}{s(\mu)}=\frac{f(\mu)-\hat f(\mu)}{\mu s(\mu)}
$$
is an entire function in the complex plane. Since the function
$g(\mu)=f(\mu)-\hat f(\mu)$ belongs to $\pwp$, it follows from
(11) that
$$
|g(\mu)|\le C_2e^{\pi|Im\mu|}.\eqno(24)
$$
It follows from (13) that
 $$
|\mu s(\mu)|\ge C_3e^{\pi|Im\mu|}\eqno(25)
$$
$(C_3>0)$ if $|Im\mu|\ge1$. This yields $|\Phi(\mu)|\le C_2/C_3$
if $|Im\mu|\ge1$.

We denote by $H$ the union of vertical segments $\{z: |Rez|=n+1/2,
|Imz|\le1\}$, where $n=N+1,N+2,\ldots$. It follows from (13) that
$|\mu s(\mu)|\ge C_4>0$ if $\mu\in H$. From the last inequality,
(24), (25) and the Maximum Principle we obtain that
$|\Phi(\mu)|\le C_5$ in the strip $|Im\mu|\le1$, hence, the
function $\Phi(\mu)$ is bounded in the whole complex plane and, by
virtue of Liouville's theorem, it is a constant. Let $|Im\mu|=1$.
Then it follows from the Paley-Wiener theorem and the Riemann
lemma that $\lim_{|\mu|\to\infty}g(\mu)=0$, hence,
$\Phi(\mu)\equiv0$.

Suppose that $q_0\ne0$. Similarly [10] we consider the function
$\tilde u(\mu)=u(\muk)$. Since $\muk=\mu+q_0/(2\mu)+O(\mu^{-3})$,
we have
$$
\cos\pi\muk=\cpm+(-\pi q_0/(2\mu)+O(\mu^{-3}))\sin\pi\mu+O(\mu^{-2})\cpm,
$$
$$
\frac{\sin\pi\muk}{\muk}=\spm+O(\mu^{-2})\cpm+O(\mu^{-3})\sin\pi\mu.
$$
This, together with (10), yields
$$
\tilde u(\mu)=\tet+\cpm+\gamma\spm+\frac{\tilde f(\mu)}{\mu},
$$
where
$$
\tilde f(\mu)=\frac{\mu f(\muk)}{\muk}+O(\mu^{-1})\cpm+O(\mu^{-2})\sin\pi\mu.
$$
It is evident that $\tilde f(\mu)\in L_2(R)$. We have proved that there exists the function
 $\tilde q(x)\in L_2(0,\pi)$ such that corresponding to the collection
$(\agt\tilde q(x))$ the characteristic determinant
$\tilde\dm=\tilde u(\mu)$ and $<\tilde q>=0$. It is readily seen
that the function $\dm=\tilde\Delta(\muq)$ is the characteristic
determinant corresponding to the collection $(\agt \tilde
q(x)+q_0)$. Since $<\tilde q+q_0>=q_0$,
$\dm=\tilde\Delta(\muq)=\tilde u(\muq)=u(\mu)$, we see that the
function $q(x)=\tilde q(x)+q_0$ has all required properties.
Theorem 2 is proved.

It is known [3, pp. 68, 80] that the eigenvalues of problem (1)+
(2) form two series
$$
\lambda_0=\mu_0^2, \quad \lambda_{n,j}=(2n+o(1))^2 \eqno(26)
$$
(in case 1) and
$$
\lambda_{n,j}=(2n-1+o(1))^2 \eqno(27)
$$
(in case 2); in both cases $j=1,2, \quad n=1,2,\ldots$. We denote
$\mu_{n,j}=\sqrt{\lambda_{n,j}}=2n-\theta+o(1)$. It follows from
[1, p. 74] that asymptotic formulas (26), (27) can be made more
precise. Namely, we have
$$
\mu_{n,j}=2n-\theta+O(n^{-1/2}). \eqno(28)
$$
It is clear that $|\mu_{n,1}-\mu_{n,2}|=O(n^{-1/2})$. If the set
of simple eigenvalues is finite, then the spectrum of problem (1)+
(2) is called asymptotically multiple. If the set of multiple
eigenvalues is finite, then the spectrum of problem (1)+(2) is
called asymptotically simple.

Many authors studied the spectrum of problems of type I [3], [4].
In particular, it is known a lot of examples of potentials $q(x)$
that ensure an asymptotically multiple spectrum. The spectrum of
any problem of type II (even if $q(x)\in L_1(0,\pi)$) is
asymptotically simple [13].

Properties of the spectrum of problems of types III and IV have been investigated
considerably less.
 If $q(x)\equiv0$, then the spectrum of any problem of type III is asymptotically multiple
and the spectrum of any problem of type IV is asymptotically
simple [2]. For problems of type III there exist examples of
potentials $q(x)\not\equiv0$ that ensure an asymptotically
multiple spectrum. It was shown in [13] that the root function
system of problem (1)+(2) with boundary conditions of type III or
IV forms a Riesz basis in $L_2(0,\pi)$ if and only if when the
spectrum is asymptotically multiple. In the author's opinion, of
great interest is the problem to give a description of the set of
potentials $q(x)$ that ensure an asymptotically multiple spectrum
for problems of type III or IV, i.e. under the condition
$\alpha\ne1/2$.

{\bf Lemma 1.} {\it Suppose that a collection $(\alpha, \gamma,
\theta, q(x))$, where $\alpha\ne1/2$, determines a characteristic
determinant $\Delta(\mu)$ represented in the form
$$
\Delta(\mu)=(-1)^{\theta+1}+\cos \pi\mu+(\gamma+\frac{\pi<q>}{2})\frac{\sin\pi\mu}{\mu}+\frac{f(\mu)}{\mu},
$$
where $f(\mu)\in PW_\pi^-$. Suppose also that the spectrum of the
Dirichlet problem with the potential $q(x)$ is simple and zero is
not an eigenvalue of this problem. Then for any $\varepsilon>0$
there exists $\delta>0$ such that for any function $\tilde
f(\mu)\in PW_\pi^-$ satisfying the condition $||f(\mu)-\tilde
f(\mu)||_{L_2(R)}<\delta$ there exists a potential $\tilde q(x)\in
L_2(0,\pi)$ such that $||q(x)-\tilde q(x)||<\varepsilon$ and the
function
$$
\tilde\Delta(\mu)=(-1)^{\theta+1}+\cos
\pi\mu+(\gamma+\frac{\pi<\tilde
q>}{2})\frac{\sin\pi\mu}{\mu}+\frac{\tilde f(\mu)}{\mu}
$$
is the characteristic determinant corresponding to the collection

\noindent$(\alpha, \gamma, \theta, \tilde q(x))$ and $<\tilde
q>=<q>$.}

Proof. Let $\mu_n$ $(n=1,2,\ldots)$ be zeros of the function $s(x,\mu)$. It is well known that
$$
\mu_n=n+O(n^{-1}).\eqno(29)
$$
It follows from (8) and (29) that for all sufficiently large $n$
$|\dot s(\pi,\mu_n)|>C_1/|\mu_n|$ $(C_1>0)$. This, together with
the simplicity of the Dirichlet problem spectrum, yields $|\dot
s(\pi,\mu_n)|>C_2/|\mu_n|$ $(C_2>0)$ for all $n$. Thus, we obtain

$$
|\mu_n\dot s(\pi,\mu_n)|>C_3>0.\eqno(30)
$$
for any $n$.

Set
$$
F(x,t)=\sum_{n=1}^\infty(u_n
\sin\mu_nx\sin\mu_nt-
\frac{2}{\pi}\sin nx\sin nt),
$$
where
$$
u_n=\frac{2c(\pi,\mu_n)}{\mu_n\dot s(\mu_n)}.$$ It is readily seen
that $F(x,t)=(G(x-t)-G(x+t))/2$, where $$
G(y)=\sum_{n=1}^\infty(u_n\cos\mu_ny-\frac{2}{\pi}\cos ny).$$ We
denote
$$
\tilde F(x,t)=\sum_{n=1}^\infty(\tilde u_n
\sin\mu_nx\sin\mu_nt-
\frac{2}{\pi}\sin nx\sin nt),
$$
where $\tilde u_n$ are some coefficients.
Trivially,
$\tilde F(x,t)=(\tilde G(x-t)-\tilde G(x+t))/2$, £¤¥ $$
\tilde G(y)=\sum_{n=1}^\infty(\tilde u_n\cos\mu_ny-\frac{2}{\pi}\cos ny).$$

We denote $g(\mu)=f(\mu)-\tilde f(\mu)$, where $\tilde f(\mu)$ is
an arbitrary function from $\pwp$ such that
$\sigma=||g(\mu)||_{L_2(R)}<1$.

Relation (29), combined with the reasons of [3, pp. 61-62],
implies that
$$
\sum_{n=1}^\infty|g(\mu_n)|^2\le C_4\sigma^2.\eqno(31)
$$

Later on we will use the following

Proposition 1. {\it If $\sum_{n=1}^\infty a_n^2<C^2$, where $a_n$ are arbitrary
nonnegative numbers, then $\sum_{n=1}^\infty\frac{a_n}{n}\le\frac{3}{2}C$.}

This proposition is a trivial corollary of the elementary
inequality $ab\le(\epsilon a^2+\epsilon^{-1}b^2)/2$
$(a,b,\epsilon>0)$ and the well-known equality $\sum_{n=1}^\infty
n^{-2}=\pi^2/6$.

We denote $R(y)=G(y)-\tilde G(y)$, $v_n=u_n-\tilde u_n$. Suppose

$$
|v_n|\le C_5|g(\mu_n)|/|\mu_n|\eqno(32)
$$
$(n=1,2,\ldots)$. Let us prove that the estimates
$$
|R(y)|\le\tilde C\sigma\eqno(33)
$$
$(0\le y\le2\pi)$ and
$$
||R'(y)||_{L_2(0,2\pi)}\le \tilde C\sigma\eqno(34)
$$
hold. At first, we will obtain estimate (33). Using (31), (32) and
proposition 1, we have
$$
|R(y)|\le\sum_{n=1}^\infty|v_n||\cos\mu_ny|\le C_6\smn|g(\mu_n)|/n\le C_7\sigma.
$$

Let us estimate the function $R'(y)$. Evidently, $R'(y)=-\smn
h_n\sin\mu_ny$, where $h_n=\mu_nv_n$. It follows from (31) and
(32) that
$$
\smn|h_n|^2\le C_8\sigma^2.\eqno(35)
$$
It is readily seen that
$$
\begin{array}{c}
||R'(y)||_\lop=||\smn h_n\sin\mu_ny||_\lop\le\\
\le||\smn h_n(\sin\mu_n y-\sin ny)||_\lop+||\smn h_n\sin
ny||_\lop.
\end{array}\eqno(36)
$$
It follows from (29), (35) and proposition 1 that
$$
\begin{array}{c}
||\smn h_n(\sin\mu_n y-\sin ny)||_\lop=\\
=2||\smn h_n\sin\frac{(\mu_n-n)y}{2}\cos\frac{(\mu_n+n)y}{2}||_\lop\le
C_9\smn\frac{|h_n|}{n}\le C_{10}\sigma.
\end{array}
$$
Evaluating the second summand in the right-hand part of (36) by
the Parseval equality, using (35), (36) and the last inequality,
we get estimate (34).

 Fix an arbitrary natural number $N_0$. Similarly, one can
prove the following assertion. If for $n>N_0$ estimate (32) holds
and for $n=1,\ldots, N_0$ the estimate
$$
|v_n|\le C_{11}\sqrt{\sigma}\eqno(32')
$$
holds, then the estimates
$$
|R(y)|\le\tilde C_1\sqrt{\sigma}\eqno(33')
$$
$(0\le y\le2\pi)$ and
$$
||R'(y)||_\lop\le\tilde C_1\sqrt{\sigma}\eqno(34')
$$
are valid.

  We denote $\Delta_+(\mu)=\dm-\tet$. Then
$$
\Delta_+(\mu_n)=\cos\pi\mu_n+(\gamma+\frac{\pi}{2}<q>)\frac{\sin\pi\mu_n}{\mu_n}+\frac{f(\mu_n)}{\mu_n}.\eqno(37)
$$
It follows from (5) that
$$
\Delta_+(\mu_n)=\alpha c(\pi,\mu_n)+(1-\alpha)s'(\pi,\mu_n).\eqno(38)
$$
We denote
$$
\tilde\Delta_+(\mu)=\cos\pi\mu+(\gamma+\frac{\pi}{2}<q>)\frac{\sin\pi\mu}{\mu}+\frac{\tilde f(\mu)}{\mu}.\eqno(39)
$$

Let us consider 3 cases.

1) $\alpha=1$. Then, it follows from (38) that
$c(\pi,\mu_n)=\Delta_+(\mu_n)$. This, together with (37), yields
$$
c(\pi,\mu_n)=\cos\pi\mu_n+(\gamma+\frac{\pi}{2}<q>)\frac{\sin\pi\mu_n}{\mu_n}+\frac{f(\mu_n)}{\mu_n}.
$$
It follows from the last equality and (39) that
$c(\pi,\mu_n)-\tilde\Delta_+(\mu_n)=g(\mu_n)/\mu_n$. Set $\tilde
c_n=\tdm$, $\tilde u_n=2\tilde c_n/(\mu_n\dot s(\pi,\mu_n))$,
then, it follows from (30) that inequality (32) holds, hence,
estimates (33) and (34) are valid.

2) $\alpha=0$. Then, it follows from (38) that
$s'(\pi,\mu_n)=\tdm$. This, combined with (4), implies that
$c(\pi,\mu_n)=1/\tdm$. Since $|c(\pi,\mu_n)|\le C_{12}$, we have
$$
|\tdm|\ge C_{12}^{-1}>0.\eqno(40)
$$
Notice, that
$$
|\Delta_+(\mu_n)-\tdm|=|g(\mu_n)|/|\mu_n|\le
C_{13}\sigma/|\mu_n|,\eqno(41)
$$
therefore, it follows from (40) and (41) that
$|\tdm|\ge(2C_{12})^{-1}$ if $\sigma$ is sufficiently small. This,
together with (40), yields
$$
|\Delta_+(\mu_n)\tdm|\ge C_{14}>0.\eqno(42)
$$
Set $\tilde c_n=1/\tdm$, $\tilde u_n=2\tilde c_n/(\mu_n\dot s(\pi,\mu_n))$,
then, it follows from (30), (41), (42) that inequality (32) holds, hence,
estimates (33) and (34) are valid.

3) $\alpha\ne1$, $\alpha\ne0$. Then, it follows from (38) and (4) that $\Delta_+(\mu_n)=
\alpha c(\pi,\mu_n)+(1-\alpha)c^{-1}(\pi,\mu_n)$, hence,
$$
\alpha c^2(\pi,\mu_n)-\Delta_+(\mu_n)c(\pi,\mu_n)+1-\alpha=0.\eqno(43)
$$
We denote $D_+(\mu)=\Delta_+^2(\mu)-4\alpha(1-\alpha)$, $\tilde D_+(\mu)=\tdm^2-4\alpha(1-\alpha)$.
Solving equation (43), we get
$$
c(\pi,\mu_n)=\frac{\Delta_+(\mu_n)+(-1)^{\delta_n}\sqrt{D_+(\mu_n)}}{2\alpha},
$$
where either $\delta_n=0$ or $\delta_n=1$. We denote
$$
\tilde c_{n,\tilde\delta_n}=\frac{\tilde\Delta_+(\mu_n)+(-1)^{\tilde\delta_n}\sqrt{\tilde D_+(\mu_n)}}{2\alpha},\eqno(44)
$$
where $\tilde\delta_n=0$ or $\tilde\delta_n=1$.
It is clear that the numbers $\tilde c_{n,\tilde\delta_n}$ are the roots of the equation $$
\alpha z^2-\tilde\Delta_+(\mu_n)z+1-\alpha=0.\eqno(45)
$$

It follows from (29), (37) and (39) that
$$
\Delta_+(\mu_n)=(-1)^n+r_n,\quad
\tilde\Delta_+(\mu_n)=(-1)^n+\tilde r_n,
$$
where
$$
r_n=\rho_n+\frac{f(\mu_n)}{\mu_n},\quad\tilde
r_n=\rho_n+\frac{\tilde f(\mu_n)}{\mu_n},
$$
where $|\rho_n|<C_{15}/n$. We denote $W_n=2(-1)^nr_n+r_n^2$,
$\tilde W_n=2(-1)^n\tilde r_n+\tilde r_n^2$. It is readily seen
that there exists a number $N_0$ such that for any $n>N_0$
$|W_n|+|\tilde W_n|<|1-2\alpha|^2/10$.

Let us consider the case $n>N_0$. Set $\tilde\delta_n=\delta_n$.
Since $D_+(\mu_n)=(1-2\alpha)^2+W_n$ and $\tilde
D_+(\mu_n)=(1-2\alpha)^2+\tilde W_n$, we see that
$$
|\sqrt{D_+(\mu_n)}+\sqrt{\tilde D_+(\mu_n)}|\ge C_{16}>0.
$$
It follows from the last inequality that
$$\begin{array}{c}
|c(\pi,\mu_n)-\tilde
c_{n,\tilde\delta_n}|=|\frac{\Delta_+(\mu_n)-\tilde\Delta_+(\mu_n)+(-1)^{\delta_n}
(\sqrt{D_+(\mu_n)}-\sqrt{\tilde D_+(\mu_n)})}{2\alpha}|\le\\\le
|\frac{\Delta_+(\mu_n)-\tilde\Delta_+(\mu_n)}{2\alpha}|+|\frac{D_+(\mu_n)-\tilde
D_+(\mu_n)}{2\alpha(\sqrt{D_+(\mu_n)}+\sqrt{\tilde
D_+(\mu_n)})}|\le\\\le
C_{17}|\Delta_+(\mu_n)-\tilde\Delta_+(\mu_n)|\le
C_{18}|g(\mu_n)|/|\mu_n|.
\end{array}
$$
Set $\tilde c_n=c_{n,\tilde\delta_n}$, $\tilde u_n=2\tilde
c_n/(\mu_n\dot s(\pi,\mu_n))$. It follows from (30) and the last
inequality that estimate (32) holds.

 Let us consider the case
$n=1,\ldots,N_0$. It follows from the equality
$(\sqrt{z_1}+\sqrt{z_2})(\sqrt{z_1}-\sqrt{z_2})=z_1-z_2$ that at
least one of the inequalities
$$
|\sqrt{z_1}+\sqrt{z_2}|\le|z_1-z_2|^{1/2},\quad
|\sqrt{z_1}-\sqrt{z_2}|\le|z_1-z_2|^{1/2}
$$ holds.
For any $n$ we choose $\tilde\delta_n$ in (44) so that the equality
$$
|(-1)^{\delta_n}\sqrt{D_+(\mu_n)}-(-1)^{\tilde\delta_n}\sqrt{\tilde
D_+(\mu_n)}|\le |D_+(\mu_n)-\tilde D_+(\mu_n)|^{1/2}
$$ is valid.
It is readily seen that
$$
|D_+(\mu_n)-\tilde D_+(\mu_n)|\le C_{19}|\Delta_+(\mu_n)-\tdm|\le
C_{20}|g(\mu_n)|/|\mu_n|.\eqno(46)
$$
Set $\tilde c_n=\tilde c_{n,\tilde\delta_n}$, $\tilde u_n=2\tilde
c_n/(\mu_n\dot s(\pi,\mu_n))$. Using (31), we obtain
$|g(\mu_n)|\le C_{21}\sigma$. Then, it follows from (30), (46) and
the last inequality that inequality (32$'$) holds, hence,
estimates (33$'$) and (34$'$) are valid.

Let us now consider the Gelfand-Levitan equation
$$
K(x,t)+\tilde F(x,t)+\int_0^xK(x,t)\tilde F(s,t)ds=0.\eqno(47)
$$
Taking into account estimates (33) and (34) (or (33$'$) and
(34$'$)) and arguing as in [10], we see that for all sufficiently
small $\sigma$ equation (47) has the unique solution $\tilde
K(x,s)$, the function $\tilde q(x)=\frac{2d\tilde K(x,x)}{dx}$
belongs to $L_2(0,\pi)$ and the inequality $||q-\tilde q||<\hat
C\sqrt{\sigma}$ holds. This implies that for all sufficiently
small $\sigma$ $||q-\tilde q||<\varepsilon$, where $\varepsilon$
is an arbitrary preassigned positive number.

We denote by $\tilde s(x,\mu)$, $\tilde c(x,\mu)$
the fundamental system of solutions to equation (1) with the potential $\tilde q(x)$ and
the initial conditions $\tilde s(0,\mu)=\tilde c'(0,\mu)=0$, $\tilde c(0,\mu)=\tilde s'(0,\mu)=1$.

As for Theorem 2, we obtain $\tilde s(\pi,\mu)\equiv s(\pi,\mu)$,
hence, the numbers $\mu_n^2$ form the Dirichlet problem spectrum
for equation (1) with the potential $\tilde q(x)$. In the same
way, we get $\tilde c(\pi,\mu_n)=\tilde c_n$. This, combined with
identity (4), yields $\tilde s'(\pi,\mu_n)=1/\tilde c_n$.

It follows from asymptotic representation (8) written for the functions $s(\pi,\mu)$ and
$\tilde s(\pi,\mu)$ that $<q>=<\tilde q>$.

Let $\hat\Delta(\mu)$ be the characteristic determinant
corresponding to the collection $(\agt\tilde q(x))$. Let us prove
that $\hat\Delta(\mu)\equiv\tilde\Delta(\mu)$. By Theorem 1, for
the function $\hat \Delta(\mu)$ we have a representation
$$
\hat\Delta(\mu)=(-1)^{\theta+1}+\cos \pi\mu+(\gamma+\frac{\pi<\tilde q>}{2})\frac{\sin\pi\mu}{\mu}+\frac{\hat f(\mu)}{\mu},
$$
where $\hat f(\mu)\in PW_\pi^-$. Taking into account (5), using
that the numbers $\tilde c_n$ are the roots of equation (45),
repeating arguments of Theorem 2, we obtain
$\hat\Delta(\mu_n)=\tilde\Delta(\mu_n)$. Then, as for Theorem 2,
we get $\hat\Delta(\mu)\equiv\tilde\Delta(\mu)$. Lemma 1 is
proved.

{\bf Lemma 2} [8]. {\it For any function $q(x)\in L_2(0,\pi)$ and
any $\varepsilon>0$ there exists a function $q_\varepsilon (x)\in
L_2(0,\pi)$ such that $||q(x)-q_\varepsilon (x)||<\varepsilon$,
the spectrum of the Dirichlet problem with the potential
$q_\varepsilon (x)$ is simple and zero is not an eigenvalue of
this problem.}

{\bf Lemma 3.} {\it Suppose the Dirichlet problem with a potential
$q(x)\in L_2(0,\pi)$ has a simple spectrum and zero is not an
eigenvalue of this problem. Then there exists $\delta>0$ such that
for any function $\tilde q(x)\in L_2(0,\pi)$ satisfying
$||q(x)-\tilde q(x)||<\delta$ the Dirichlet problem with the
potential $\tilde q(x)$ also has a simple spectrum and zero is not
an eigenvalue of this problem.}

Proof. Let $\lambda_n=\mu_n^2$ $(Re\mu_n\ge0)$ be the eigenvalues
of the Dirichlet problem with the potential $q(x)$. It is clear
that there exists $c_1>0$ such that $|\mu_n|>c_1$ and
$|\mu_i-\mu_j|>c_1$ if $i\ne j$. Let $s(x,\mu)$ be the solution to
equation (1) defined above,
$\Gamma=\bigcup_{n=1}^\infty\Gamma(\mu_n,c_1/3)$. Then [14, p. 87]
outside $\Gamma$ we have the estimate
$$
|s(\pi,\mu)|\ge c_2/(|\mu|+1).\eqno(48)
$$
Let $\tilde q(x)$ be a function such that $||\tilde q(x)||\le
2||q||+1$.
 Let $\tilde s(x,\mu)$ be the solution to the equation $u''-\tilde q(x)u+\mu^2u=0$
satisfying the same initial conditions as the function $s(x,\mu)$.
Then, it follows from [15, p. 46] and the trivial inequality
$\int_0^\pi|\tilde q(x)|dx\le\sqrt{\pi}||\tilde q||$ that for all
$\mu$ such that
$$
|\mu|>2\sqrt{\pi}||\tilde q||+1\eqno(49)
$$the estimate
$$
|\tilde
s(\pi,\mu)-\spm|\le\frac{c_3e^{|Im\mu|\pi}}{|\mu|^2}\eqno(50)
$$
is valid (a constant $c_3$ does not depend on $\tilde q(x)$).

Let $\tilde\mu$ be a root of the function $\tilde s(\pi,\mu)$ .
Suppose $|Im\tilde\mu|>1$, then $|\sin\pi\tilde\mu|\ge
c_4e^{|Im\tilde\mu|\pi}$ $(c_4>0)$. Suppose also that $\tilde\mu$
satisfies (49), then it follows from the last inequality and (50)
that $|\tilde\mu|<c_3/c_4$, hence, $|Im\tilde\mu|<c_3/c_4$. This
implies that outside of the strip $|Im\mu|\le c_3/c_4$ the
function $\tilde s(\pi,\mu)$ has no roots.

 Let us study properties of the function $\tilde s(\pi,\mu)$ inside the strip
$|Im\mu|\le c_3/c_4+c_1/3$. It follows from (6) and [7]
$$
|s(\pi,\mu)-\tilde s(\pi,\mu)|\le c_5||q-\tilde q||/(|\mu|+1).\eqno(51)
$$
It follows from (48) and (51) that under the condition
$$
||q-\tilde q||\le c_2/(2c_5)\eqno(52)
$$
the function $\tilde s(\pi,\mu)$ has no roots outside of $\Gamma$,
therefore, all the roots of the function $\tilde s(\pi,\mu)$ lie
inside $\Gamma$. Combining (48), (51) and the Rouch$\acute e$
theorem, we obtain that under condition (52) the function $\tilde
s(\pi,\mu)$ has a unique root inside each disk $\Gamma(\mu_n,
c_1/3)$. This completes the proof of lemma 3.

It follows from equality (5), asymptotic formula (28) and the
Hadamard theorem [16, p. 259] that
$$
\Delta(\mu)=\frac{\pi^2}{2}(\mu_0^2-\mu^2)\prod_{n=1}^\infty\frac{(\mu_{n,1}^2-\mu^2)(\mu_{n,2}^2-\mu^2)}{16n^4}\eqno(53)
$$
if $\theta=0$ and
$$
\Delta(\mu)=2\prod_{n=1}^\infty\frac{(\mu_{n,1}^2-\mu^2)(\mu_{n,2}^2-\mu^2)}{(2n-1)^4}\eqno(54)
$$
if $\theta=1$.

{\bf Lemma 4} [6]. {\it Suppose the numbers $\mu_{n,j}$ satisfy
asymptotic relations
$$
\mu_{n,j}=2n-\theta+\frac{B_j}{2n-\theta}+\frac{\delta_{n,j}}{n},
$$
where $\sum_{n=1}^\infty|\delta_{n,j}|^2<\infty$ $(j=1,2)$, then the function
$\Delta(\mu)$ can be represented in the form
$$
\Delta(\mu)=(-1)^{\theta+1}+\cos \pi\mu+\frac{\pi}{2}(B_1+B_2)\frac{\sin\pi\mu}{\mu}+\frac{f(\mu)}{\mu},
$$
where $f(\mu)\in PW_\pi^-$.}

Consider two sequences
$$
\mu_{n,j}=2n-\theta+\frac{V_1}{2n-\theta}+\frac{V_2}{(2n-\theta)^2}+\frac{\delta_{n,j}}{n^2},
$$
where $V_1$ and $V_2$ are arbitrary complex numbers and
$\sum_{n=1}^\infty|\delta_{n,j}|^2<\infty$ $(j=1,2)$. Let $\mu_0$
be an arbitrary complex number and let a function $\Delta(\mu)$ be
determined by (53) or (54). We denote $\tilde\mu_0=\mu_0$,
$\tilde\mu_{n,j}=\mu_{n,j}$ if $n=\overline{1, N}$ and
$$\tilde\mu_{n,1}=\tilde\mu_{n,2}=2n-\theta+\frac{V_1}{2n-\theta}+\frac{V_2}{(2n-\theta)^2}$$
if $n=N+1, N+2,\ldots$, where $N$ is an arbitrary number. Set
$$
\Delta_N(\mu)=\frac{\pi^2}{2}(\tilde\mu_0^2-\mu^2)\prod_{n=1}^\infty\frac{(\tilde\mu_{n,1}^2-\mu^2)(\tilde\mu_{n,2}^2-\mu^2)}{16n^4}\eqno(55)
$$
if $\theta=0$ and
$$
\Delta_N(\mu)=2\prod_{n=1}^\infty\frac{(\tilde\mu_{n,1}^2-\mu^2)(\tilde\mu_{n,2}^2-\mu^2)}{(2n-1)^4}\eqno(56)
$$
if $\theta=1$.

{\bf Lemma 5.} {\it In both cases
$$
\lim_{N\to\infty}||\mu(\Delta(\mu)-\Delta_N(\mu))||_{L_2(R)}=0.
$$}

Proof. By lemma 4, the function $\dm$ can be represented in the
form
$$
\dm=\tet +\cpm+\pi V_1\spm+\frac{f(\mu)}{\mu}, \eqno(57)
$$
where $f(\mu)\in\pwp$, and the function $\dnm$ can be represented in the form
$$
\dnm=\tet+\cpm+\pi V_1\spm+\frac{f_N(\mu)}{\mu},\eqno(58)
$$
where $f_N(\mu)\in\pwp$. We denote $F_N(\mu)=\mu(\dm-\dnm)$. It follows from (57) and (58)
that $F_N(\mu)=f(\mu)-f_N(\mu)$. It is readily seen that
$$
|Im\mun|\le M,\quad|Im\tilde \mun|\le M,\quad|\mun^2-\tilde\mun^2|<C_0.\eqno(59)
$$
We denote by $l$ the line $Imz=(M+2)(C_0+1)$. Let $\mu\in l$, then $\dm\ne0$
and, hence, we have
$$
F_N(\mu)=\mu\dm\left(1-\frac{\dnm}{\dm}\right)=\mu\dm\left(1-\prod_{n=N+1}^\infty\prod_{j=1}^2(1+\frac{\tilde\mun^2-\mun^2}{\mun^2-\mu^2})\right).\eqno(60)
$$
It follows from (57) and (11) that
$$
|\dm|<c_1.\eqno(61)$$
We denote $\aln=\frac{\tilde\mun^2-\mun^2}{\mun^2-\mu^2}$ $(j=1,2)$. It follows from (59)
that
$$
|\aln|<1/4.\eqno(62)
$$
We denote $\wnj=\sum_{n=N+1}^\infty\ln(1+\aln)$. Here we take the
branch of $\ln(1+z)$ which vanishes for $z=0$. Let us estimate the
functions $\aln$ and $\wnj$. First of all, notice, that they are
even, therefore, one can consider only the case $Re\mu\ge0$.
Obviously,
$$
|\tilde \mun^2-\mun^2|\le c_2|\delta_n|/n.\eqno(63)
$$
It is not hard to prove that for all $n>N_0$, where $N_0$ is a sufficiently large number, we have

$$
|\mun+\mu|>|\mu|,\quad|\mun-\mu|>\ntm/2.\eqno(64)
$$
 It follows from (63) and (64) that if $n>N_0$, then
$$
|\aln|\le\frac{2c_2|\delta_n|}{n|\mu|\ntm}.\eqno(65)
$$
Further we count that $N>N_0$. Using (62), (65) and the elementary
inequality $|\ln(1+z)|<2|z|$, which holds if $|z|<1/4$, we obtain that
$$
\begin{array}{c}
|\wnj|\le\sum_{n=N+1}^\infty|\ln(1+\aln)|\le\\
\le2\sum_{n=N+1}^\infty|\aln|\le\frac{4c_2}{|\mu|}\sum_{n=N+1}^\infty\frac{|\delta_n|}{n\ntm}.
\end{array}
\eqno(66)
$$
It follows from (66) and the trivial inequality
$|\delta_n|/n\le(|\delta_n|^2+n^{-2})/2$ that for all sufficiently
large $N$ $|\wnj|<1/8$, hence, $\sum_{j=1}^2|\wnj|<1/4$. This,
together with the elementary inequality $|1-e^z|<2|z|$, which
holds if $|z|<1/4$, yields
$$
|1-\exp(\sum_{j=1}^2\wnj)|<2\sum_{j=1}^2|\wnj|.
$$
It follows from the last inequality, (60), (61) and (66) that
$$
|F_N(\mu)|\le c_3|\mu|\sum_{j=1}^2|\wnj|\le c_4\sum_{n=N+1}^\infty\frac{|\delta_n|}{\ntm}.\eqno(67)
$$
To estimate the sum in the right-hand part of (67) we need the elementary inequality
$$
\sum_{n=1}^\infty|n-z|^{-2}<\tilde C,\eqno(68)
$$
where $\tilde C$ does not depend of $z\in l$. Fix an arbitrary
$\varepsilon>0$. We choose $\varepsilon_0>0$ small enough for the
inequality $\varepsilon_0\tilde C<\varepsilon/10$ to be valid.
Using (68), we obtain
$$
\begin{array}{c}
\sni\frac{|\delta_n|}{n\ntm}\le\frac{1}{|\mu|}\sni\frac{|\delta_n|(\ntm+2n-\theta)}{n\ntm}|\le\\
\le\frac{1}{|\mu|}[\sni|\delta_n|/n+\sni\frac{|\delta_n|(2n-\theta)}{n\ntm}|]\le\\
\le\frac{1}{|\mu|}[\sni(|\delta_n|^2+n^{-2})/2+2\sni\frac{|\delta_n|}{\ntm}|]\le\\
\le\frac{1}{|\mu|}[(1/2+1/\varepsilon_0)\sni|\delta_n|^2+\\+1/2\sni n^{-2}+\varepsilon_0\sni\frac{1}{\ntm}]\le\\
\le\frac{1}{|\mu|}[(1/2+1/\varepsilon_0)\sni|\delta_n|^2+1/2\sni n^{-2}+\varepsilon/10].
\end{array}
$$
Evidently, that for all sufficiently large $N$ the right-hand part
of the last inequality does not exceed $\varepsilon/|\mu|$. This,
combined with (67), implies that for $\mu\in l$ $|F_N(\mu)|\le
C_N/|\mu|$, where $C_N\to0$ as $N\to\infty$. This yields
$$
\lim_{N\to\infty}||F_N(\mu)||_{L_2(l)}=0.
$$
It follows from the last relation and [9, p. 115] that
$$
\lim_{N\to\infty}||F_N(\mu)||_{L_2(R)}=0.
$$
Lemma 5 is proved.

{\bf Lemma 6.} {\it Let $\alpha\ne1/2$. Let $p$ be an arbitrary
odd number, $p=2l+1$, where $l=0,1,\ldots$. Then there exist
numbers $h_i, g_i$ $(i=\overline{0,p-1})$ such that for any
function $q(x)\in W_1^p[0,\pi]$
 satisfying the conditions $q^{(i)}(0)=h_i$, $q^{(i)}(\pi)=g_i$
 $(i=\overline{0,p-1})$ we have for the numbers $\mu_{n,j}$ the asymptotic representation
$$
\mu_{n,j}=2n-\theta+\sum_{m=1}^{l+1}\frac{V_m}{(2n-\theta)^m}+o(n^{-l-1})\eqno(69)
$$
$(j=1,2)$ , moreover, $V_1=\pi^{-1}(\gamma+\pi<q>/2)$.}

Proof. Let $q(x)\in W_1^p[0,\pi]$. It was shown in [3, p. 69] that
the characteristic equation $\Delta(\mu)=0$ can be reduced to the
form
$$
\{\exp[i\pi\mu+\int_0^\pi\sigma(\mu,t)dt]+Cw(\mu,0)/G(\mu)\}^2=H(\mu)/(w(\mu,0)w(\mu,\pi)), \eqno(70)
$$
where
$$
C=-(A_{13}+A_{24})/2, \eqno(71)
$$
$$
G(\mu)=-i\mu(A_{14}+A_{23})+A_{14}\sigma(-\mu,0)-A_{23}\sigma(\mu,\pi)+A_{34},\eqno(72)
$$
$$
w(\mu,x)=2i\mu+\sigma(\mu,x)-\sigma(-\mu,x),\eqno(73)
$$
$$
H(\mu)=C^2w(\mu,0)w(\mu,\pi)+G(\mu)G(-\mu),\eqno(74)
$$
$$
\sigma(\mu,x)=\sum_{k=1}^p\frac{\sigma_k(x)}{(2i\mu)^k}+\frac{\sigma_p(\mu,x)}{(2i\mu)^p},\eqno(75)
$$
$$
\sigma_p(\mu,x)=\int_0^x\sigma_{p+1}(x-\xi)e^{-2i\mu\xi}d\xi+o(\mu^{-1})\eqno(76)
$$
[3, pp. 60, 61, 69], and the functions $\sigma_k(x)$ are
determined by the recursion relations
$$
\sigma_1(x)=q(x), \sigma_2(x)=-q'(x),
\sigma_{k+1}=-\sigma_k'(x)-\sum_{j=1}^{k-1}\sigma_{k-j}(x)\sigma_j(x).\eqno(77)
$$
It follows from (77) that
$$
\sigma_{k+1}(x)=(-1)^kq^{(k)}(x)+S_{k-2}(x),\eqno(78)
$$
where $S_{k-2}(x)$ is a polynomial of $q(x), q'(x),\ldots, q^{(k-2)}(x)$. It follows from (73)
that
$$
w(\mu,x)=2i\mu+O(\mu^{-1}).\eqno(79)
$$
Using relations (71)-(74) and performing some simple though
awkward manipulations, we obtain
$$
\begin{array}{c}
H(\mu)=\frac{1}{4} (A_{13}+A_{24})^2\times\\\times[2i\mu
+\sigma(\mu,0)-\sigma(-\mu,0)]
[2i\mu+\sigma(\mu,\pi)-\sigma(-\mu,\pi)]+\\
+[-i\mu(A_{14}+A_{23})+A_{14}\sigma(-\mu,0)-A_{23}\sigma(\mu,\pi)+A_{34}]
\times\\
\times[i\mu(A_{14}+A_{23})
+A_{14}\sigma(\mu,0)
-A_{23}\sigma(-\mu,\pi)
+A_
{34}
]=\\
=\frac{i\mu}{2}(A_{13}+A_{24})^2[\sigma(\mu,\pi)
-\sigma(-\mu,\pi)+\sigma(\mu,0)-\sigma(-\mu,0)]+\\
+i\mu(A_{14}+A_{23})[A_{14}\sigma(-\mu,0)
-A_{23}\sigma(\mu,\pi)-\\-A_{14}\sigma(\mu,0)
+A_{23}\sigma(-\mu,\pi)]+\\
+A_{34}^2+\frac{1}{4}(A_{13}+A_{24})^2
[\sigma(\mu,0)-\sigma(-\mu,0)]
[\sigma(\mu,\pi)-\sigma(-\mu,\pi)]+\\+[
A_{14}\sigma(\mu,0)-A_{23}\sigma(-\mu,\pi)][A_{14}\sigma(-\mu,0)-A_{23}\sigma(\mu,\pi)].
\end{array}\eqno(80)
$$
We denote $\sum_{k=1}^p\sigma_k(x)/(2i\mu)^k=\hat\sigma_p(\mu,x)$.

It follows from (75), (76) and (80) that
$$
\begin{array}{c}
H(\mu)=\frac{i\mu}{2}(A_{13}+A_{24})^2[\hat\sigma_p(\mu,\pi)
-\hat\sigma_p(-\mu,\pi)+\hat\sigma_p(\mu,0)-\hat\sigma_p(-\mu,0)]+\\
+i\mu(A_{14}+A_{23})[A_{14}\hat\sigma_p(-\mu,0)
-A_{23}\hat\sigma_p(\mu,\pi)-\\-A_{14}\hat\sigma_p(\mu,0)
+A_{23}\hat\sigma_p(-\mu,\pi)]+\\
+A_{34}^2+\frac{1}{4}(A_{13}+A_{24})^2
[\hat\sigma_p(\mu,0)-\hat\sigma_p(-\mu,0)]
[\hat\sigma_p(\mu,\pi)-\hat\sigma_p(-\mu,\pi)]+\\+[
A_{14}\hat\sigma_p(\mu,0)-A_{23}\hat\sigma_p(-\mu,\pi)][A_{14}\hat\sigma_p(-\mu,0)-A_{23}\hat\sigma_p(\mu,\pi)]+o(\mu^{1-p}).
\end{array}\eqno(81)
$$

We denote the sum of the first three summands in the right-hand
part of the last equality by $H_1(\mu)$ and we denote the sum of
the fourth and the fifth summands by $H_2(\mu)$. It follows from
(3), (75), (77) that
$$
\begin{array}{c}
H_1(\mu)=\frac{1}{2}(A_{13}+A_{24})^2(\sigma_1(\pi)+\sigma_1(0))-\\
-(A_{14}+A_{23})(A_{14}\sigma_1(0)+A_{23}\sigma_1(\pi))+A_{34}^2+\\
+\sum_{k=2}^p\frac{1-(-1)^k}{2(2i\mu)^{k-1}}[(A_{13}+A_{24})^2(\sigma_k(\pi)+\sigma_k(0))/2-\\
-(A_{14}+A_{23})(A_{14}\sigma_k(0)+A_{23}\sigma_k(\pi)]=\\
=\{\frac{1}{2}(A_{14}+A_{23})(A_{14}-A_{23})(q(\pi)-q(0))+A_{34}^2\}+\\
+\sum_{m=1}^l(2i\mu)^{-2m}(A_{14}+A_{23})(A_{14}-A_{23})(\sigma_{2m+1}(\pi)-\sigma_{2m+1}(0))/2
\end{array}\eqno(82)
$$
(if $l=0$ the last summand in the right-hand part of (82) is missing).

It is readily seen that $H_2(\mu)$ is an even function. This,
together with (75), yields
$$\begin{array}{c}
H_2(\mu)=\sum_{m=1}^{2l+1}\mu^{-2m}\sum_{i+j=2m,\ 1\le i,j\le p}
(\alpha_{ij}\sigma_i(0)\sigma_j(\pi)+\\+
\beta_{ij}\sigma_i(0)\sigma_j(\pi)+\gamma_{ij}\sigma_i(\pi)\sigma_j(\pi)),
\end{array}\eqno(83)
$$
where $\alpha_{ij}$, $\beta_{ij}$, $\gamma_{ij}$ are some
coefficients. We denote the expression in braces in the right-hand
part of (82) by $H_0$.

Let us consider the case $l=0$. It follows from (81)-(83) that
$H(\mu)=H_0+H_2(\mu)+o(1)=H_0+o(1)$. This, combined with (79), implies that
$$
H(\mu)/(w(\mu,0)w(\mu,\pi))=O(\mu^{-2}) \eqno(84)
$$
and if $H_0=0$, then
$$
H(\mu)/(w(\mu,0)w(\mu,\pi))=o(\mu^{-2}). \eqno(85)
$$
Consider the left-hand part of (70). It follows from (75-77) that
$$
\int_0^\pi\sigma(\mu,t)dt=(2i\mu)^{-1}\int_0^\pi q(x)dx+o(\mu^{-1}).\eqno(86)
$$
Using (3), (71), (72), (79) and the equality
$$
(1+c/\mu+O(\mu^{-2}))^{-1}=1-c/\mu+O(\mu^{-2})
$$
($c$ is an arbitrary number), which can easily be checked, we
obtain
$$
\begin{array}{c}
\frac{Cw(\mu,0)}{G(\mu)}=\frac{-(A_{13}+A_{24})(2i\mu+O(\mu^{-1}))}{2[-i\mu(A_{14}+A_{23})+A_{14}\sigma(-\mu,0)-A_{23}\sigma(\mu,\pi)+A_{34}]}=\\
=(-1)^{\theta+1}+\frac{A_{34}}{i\mu(A_{13}+A_{24})}+O(\mu^{-2})=(-1)^{\theta+1}(1-\frac{\gamma}{i\mu}+O(\mu^{-2})).
\end{array}\eqno(87)
$$
Rough asymptotic relation (28) yields that equation (70) has a root
$\mu_k=k+\varepsilon_k$, where $\varepsilon_k=O(k^{-1/2})$, and even $k$
correspond to case 1 $(\theta=0)$, and odd $k$ correspond to case 2 $(\theta=1)$.
This, combined with (84), (86) and (87), implies that $\varepsilon_k$ is a root
of the equation
$$
i\pi\varepsilon_k+(2i\mu_k)^{-1}\int_0^\pi
q(x)dx=-\gamma/(i\mu_k)+O(\mu_k^{-1}).
$$
This yields that $\varepsilon_k=O(k^{-1})$, and if $H_0=0$, then,
using (85), we get
$$
\varepsilon_k=(\pi\mu_k)^{-1}(1/2\int_0^\pi q(x)dx+\gamma)+o(k^{-1})=(\pi k)^{-1}(\pi<q>/2+\gamma)+o(k^{-1}).
$$
It follows from the last equality that if $H_0=0$, then for the numbers $\mu_{n,j}$
we have the asymptotic representation
$$
\mu_{n,j}=2n-\theta+\frac{\pi^{-1}(\pi<q>/2+\gamma)}{2n-\theta}+o(n^{-1}).\eqno(88)
$$
Notice, that $H_0=0$ if
$$
q(\pi)-q(0)=\frac{2A_{34}^2}{(A_{14}+A_{23})(A_{23}-A_{14})}.\eqno(89)
$$

Let us consider the general case when $l$ is an arbitrary natural
number. Let condition (89) hold. Combining (77), (78), (81), (82)
and (83), we obtain
$$
\begin{array}{c}
H(\mu)=\sum_{m=1}^l\mu^{-2m}[(2i)^{2m}(A_{14}+A_{23})(A_{14}-A_{23})(q^{(2m)}(\pi)-q^{(2m)}(0))/2+\\
+P_{2m-2}(0,\pi)]+o(\mu^{-2l}),
\end{array}
$$
where $P_{2m-2}(0,\pi)$ is a polynomial of $q(0), \ldots, q^{(2m-2)}(0), q(\pi), \ldots , q^{(2m-2)}(\pi)$.
Setting consecutively for $m=1,2,\ldots, l$
$$
q^{(2m)}(\pi)-q^{(2m)}(0)=\frac{2(2i)^{2m}P_{2m-2}(0,\pi)}{(A_{14}+A_{23})(A_{23}-A_{14})},\eqno(90)
$$
we obtain that under the conditions (89), (90) $H(\mu)=o(\mu^{-2l})$. This, together with (79), yields that
equation (70) takes the form
$$
\{\exp[i\pi\mu+\int_0^\pi\sigma(\mu,t)dt]+Cw(\mu,0)/G(\mu)\}^2=o(\mu^{-2l-2}).\eqno(91)
$$
Let us consider the left-hand part of (91). It follows from (75) that
$$
\int_0^\pi\sigma(\mu,t)dt=\sum_{j=1}^{2l+1}c_j(2i\mu)^{-j}+o(\mu^{-2l-1}),\eqno(92)
$$ where $c_j$ are some coefficients.
It follows from (73), (75) that $w(\mu,0)=2i\mu(1+\varphi(\mu))$, where
$$
\varphi(\mu)=\sum_{j=1}^l\frac{\sigma_{2j+1}(0)}{(2i\mu)^{2j+1}}+o(\mu^{-2l-2}).\eqno(93)
$$
Combining (71), (72) and (92), we obtain
$$
\begin{array}{c}
Cw(\mu,0)/G(\mu)=-i\mu(A_{13}+A_{24})(1+\varphi(\mu)/G(\mu)=\\
=\frac{(A_{13}+A_{24})(1+\varphi(\mu))}{A_{14}+A_{23}-(A_{34}+A_{14}\sigma(-\mu,0)-A_{23}\sigma(\mu,\pi))/(i\mu)}=
\frac{(-1)^{\theta+1}(1+\varphi(\mu))}{1+\psi(\mu)},
\end{array}
\eqno(94)
$$
where
$$
\psi(\mu)=\frac{-A_{34}}{i\mu(A_{14}+A_{23})}(A_{34}+A_{14}\sigma(\mu,0)-A_{23}\sigma(\mu,\pi)).\eqno(95)
$$
It was shown above that equation (91) has a root
$\mu_k=k+\varepsilon_k$, where $\varepsilon_k=O(k^{-1})$, and even
$k$ correspond to case 1 $(\theta=0)$ and odd $k$ correspond to
case 2 $(\theta=1)$. It follows from (91), (95) that
$$
\exp[i\pi\mu_k+\int_0^\pi\sigma(\mu_k,t)dt]=\frac{(-1)^{\theta}(1+\varphi(\mu_k))}{1+\psi(\mu_k)}+
o(\mu_k^{-l-1}).\eqno(96)
$$
Using (95) and (96), we get
$$
i\pi\varepsilon_k+\int_0^\pi\sigma(\mu_k,t)dt=\ln[1+\varphi(\mu_k)+o(\mu^{-l-1})]-\ln[1+\psi(\mu_k)].
$$
It follows from (75), (94), (95) and the last equality that
$$
\varepsilon_k=\left.F(\omega)\right|_{\omega=(k+\varepsilon_k)^{-1}}+o(k^{-l-1}),
$$
where $F(\omega)=\sum_{j=1}^{l+1}f_j\omega^j$, where $f_j$ are
some coefficients. Arguing as in [3, pp. 72-75], we obtain
relation (69).

Lemma 6 is proved.

{\bf Lemma 7.} {\it For any function $f(x)\in L_2(K)$, where
$K=[a,b]$, any $\varepsilon>0$ and any numbers $h_i$, $g_i$
$(i=0,\ldots,m-1)$, where $m$ is an arbitrary natural number,
there exists a function $\tilde f(x)\in C^\infty(K)$ such that
$\tilde f^{(i)}(a)=h_i$, $\tilde f^{(i)}(b)=g_i$
$(i=0,\ldots,m-1)$ and $||f(x)-\tilde
f(x)||_{L_2(K)}<\varepsilon$.}

Proof. Evidently, there exists a trigonometric polynomial $T(x)$
such that
$$
||f(x)-T(x)||_{L_2(K)}<\varepsilon/4, \eqno(97)
$$
and there exist numbers $a_0$ and $b_0$ $(a<a_0<b_0<b)$ such that
$$
||T(x)||_{L_2(K-K_0)}<\varepsilon/4, \eqno(98)
$$
where $K_0=[a_0,b_0]$. Let $\eta(x)$ be the cut-off function:
$\eta(x)=1$ if $x\in K_0$, $\eta(x)=0$ if $x\notin K$,
$0\le\eta(x)\le1$, $\eta(x)\in C^\infty(-\infty,\infty)$. We
denote $P_1(x)=\sum_{n=1}^{m-1}\frac{h_n}{n!}(x-a)^n$,
$P_2(x)=\sum_{n=1}^{m-1}\frac{g_n}{n!}(x-b)^n$. It is obvious that
$P_1^{(i)}(a)=h_i$, $P_2^{(i)}(b)=g_i$ $(i=0,\ldots,m-1)$.
Evidently, there exist segments $K_1=[a,a_1]$ and $K_2=[b_1,b]$
$(a<a_1<b_1<b)$ such that
$$
||P_j(x)||_{L_2(K_j)}<\varepsilon/4\eqno(99)
$$
$(j=1,2)$. Let us define the cut-off functions $\eta_j(x)$: $\eta_1(x)=1$ if
$x\le a$, $\eta_1(x)=0$ if $x\ge a_1$, $\eta_2(x)=1$ if $x\ge b$, $\eta_2(x)=0$
if $x\le b_1$, $0\le\eta_j(x)\le1$, $\eta_j(x)\in C^\infty(-\infty,\infty)$.

Set $\tilde f(x)=T(x)\eta(x)+P_1(x)\eta_1(x)+P_2(x)\eta_2(x)$. It is readily seen that
$\tilde f^{(i)}(a)=h_i$, $\tilde f^{(i)}(b)=g_i$ $(i=0,\ldots,m-1)$. It follows from (97)-(99)
that
$$\begin{array}{c}
||f(x)-\tilde f(x)||_{L_2(K)}\le||f(x)-T(x)||_{L_2(K)}+||(1-\eta(x))T(x)||_{L_2(K-K_0)}+\\
+\sum_{j=1}^2||P_j(x)\eta_j(x)||_{L_2(K_j)}<\varepsilon.
\end{array}
$$
Lemma 7 is proved.

The main result of the present paper is the following

{\bf Theorem 3.} {\it If $\alpha\ne1/2$, then for any
$\varepsilon>0$ there exists a function $\tilde q(x)\in
L_2(0,\pi)$ such that $||q(x)-\tilde q(x)||<\varepsilon$ and
problem (1)+(2) with the potential $\tilde q(x)$ has an
asymptotically multiple spectrum.}

Proof. Fix an arbitrary $\varepsilon>0$. By lemma 2, there exists
a function $q^{[1]}(x)\in L_2(0,\pi)$ such that
$||q-q^{[1]}||<\varepsilon/10$
 and the Dirichlet problem with the potential $q^{[1]}(x)$ has a simple spectrum, moreover, zero
is not an eigenvalue of this problem. By virtue of lemma 3, there
exists $\delta$ such that $0<\delta<\varepsilon/10$ and for any
function $\hat q(x)\in L_2(0,\pi)$ satisfying the condition
$||q^{[1]}-\hat q||<\delta$ the Dirichlet problem with the
potential $\hat q(x)$ also has a simple spectrum, moreover, zero
is not an eigenvalue of this problem.

According to lemmas 6 and 7, there exists a function
$q^{[2]}(x)\in C^\infty[0,\pi]$ such that $||q^{[1]}-\hat q||$ and
for the eigenvalues $\lambda_0=\mu_0^2$ (if $\theta=0$),
$\lambda_{n,j}=\mun^2$
 of problem (1)+(2) with the potential $q^{[2]}(x)$
we have asymptotic relation (69), where $l=2$ and
 $V_1=\pi^{-1}(\gamma+\pi<q^{[2]}>/2)$. For the characteristic determinant
 $\dm$ of the mentioned problem relations (53), (54)and (57) are valid. Let the functions
$\dnm$ and $f_N(\mu)$ be determined by formulas (55), (56), (58).
By lemma 5, it follows that
$$
\lim_{N\to\infty}||f(\mu)-f_N(\mu)||_{L_2(R)}=0.\eqno(100)
$$
Since the Dirichlet problem with the potential $q^{[2]}(x)$ has a
simple spectrum and zero is not an eigenvalue of this problem, we
see that it follows from
 (100) and lemma 1 that for any $N>N_0$, where $N_0$ is a sufficiently large
number, there exists a potential $q_N(x)\in L_2(0,\pi)$ such that
the function $\dnm$ is the characteristic determinant
corresponding to the collection $(\agt q_N(x))$, moreover,
$\lim_{N\to\infty}||q^{[2]}-q_N||=0$. Since any potential $q_N(x)$
$(N>N_0)$ ensures an asymptotically multiple spectrum of
corresponding problem (1)+(2), we see that from the last equality
it follows that theorem 3 is valid.

In the case $\alpha=1/2$, $\gamma=0$ the analogous proposition was
obtained in [17].

\medskip
\medskip

\centerline{Acknowledgments}
\medskip

The author thanks Academician V.A. Il'in for valuable discussions
of the results.

\medskip
\medskip

\centerline {\bf References}

\medskip
\medskip
[1] M.A. Naimark, Linear Differential Operators (Moscow, 1969).

[2] P. Lang, J. Locker, Spectral theory of two-point differential
operators determined by $-D^2$, J. Math. Anal. Appl. {\bf 146},
148-191 (1990).

[3] V.A. Marchenko, Sturm-Liouville Operators and Their
Applications (Kiev, 1977).

[4] F. Gesztesy and R. Weikard, Elliptic algebro-geometric
solutions of the KdV and AKNS hierarchies -- an analytic approach,
Bull. Amer. Math. Soc. {\bf 35}, 271-317 (1998).

[5] O.A. Plaksina, Characterization of the spectrum of some
boundary value problems for the Sturm-Liouville operator, Uspekhi
Mat. Nauk {\bf 36}, 199-200 (1981).

[6] I.M. Guseinov, I.M. Nabiev, Solutions of a class of inverse
Sturm-Liouville boundary value problems, Matem. Sb. {\bf 186},
35-48 (1995).

[7] L.A. Pastur, V.A. Tkachenko, Spectral theory of a class of
one-dimensional Schrodinger operators with limit-periodic
potentials, Trudy Moscow. Mat. Obshch. {\bf 51}, 114-168 (1988).

[8] J.-J. Sansuc and V. Tkachenko, Spectral parametrization of
non-selfadjoint Hill's operators, J. Diff. Eq. {\bf 125}, 366-384
(1996).

[9] S.M. Nikolskii, Approximation of Functions of Many Variables
and Embedding Theorems (Moscow, 1977).

[10] V.A. Tkachenko, Discriminants and generic spectra of
nonselfadjoint Hill's operators, Adv. Sov. Math. {\bf 19}, 41-71
(1994).

[11] A.M. Sedletskii, The stability of the completeness and of the
minimality in $L_2$ of a system of exponential functions, Mat.
Zametki {\bf 15}, 213-219 (1974).

[12] A.M. Sedletskii, Convergence of anharmonic Fourier series in
systems of exponentials, cosines and sines, Dokl. Akad. Nauk SSSR
{\bf 301}, 1053-1056 (1988).

[13] A.S. Makin, On spectral decompositions corresponding to
non-self-adjoint Sturm-Liouville operators, Dokl. Akad. Nauk {\bf
406}, 21-24 (2006).

[14] V.Ya. Levin, I.V. Ostrovskii, Small perturbations of the set
of roots of sine-type functions, Izv. Akad. Nauk SSSR, Ser. Mat.
{\bf 43}, 87-110 (1979).

[15] A.G. Kostuchenko, I.S. Sargsyan, Distribution of Eigenvalues
(Moscow, 1979).

[16] E. Titchmarsh, The Theory of Functions (Moscow, 1980).

[17] V.A. Tkachenko, Spectral analysis of non-selfadjoint Hill
operator, Dokl. Akad. Nauk SSSR {\bf 322}, 248-252 (1992).

\medskip
\medskip

Moscow State Academy of Instrument-Making and Informatics,
Stromynka 20, Moscow, 107996, Russia
\medskip

E-mail address: alexmakin@yandex.ru

}
\end{document}